\documentstyle[leqno]{article}
\makeatletter
\def\@begintheorem#1#2{\sl \trivlist \item[\hskip \labelsep{\bf #1\ #2.}]}
\def\@opargbegintheorem#1#2#3{\sl \trivlist
      \item[\hskip \labelsep{\bf #1\ #2\ (#3).}]}
\def\section{\@startsection {section}{0}{\z@}{-3.5ex plus -1ex minus 
    -.2ex}{2.3ex plus .2ex}{\bf}}
\makeatother
\newtheorem{th}{Theorem}[section]
\newtheorem{cor}{Corollary}[section]
\newtheorem{ex}{Example}[section]
\newtheorem{lm}{Lemma}[section]
\newtheorem{defn}{Definition}[section]
\newtheorem{prop}{Proposition}[section]
\newtheorem{rem}{Remark}[section]
\renewcommand{\a}{\alpha}
\renewcommand{\b}{\beta}
\newcommand{\qed}{\hbox{\rule{6pt}{6pt}}}
\begin{document}
\title{\normalsize {\bf On spectra of noises associated with Harris flows}}
\author{\normalsize  Jon Warren and Shinzo Watanabe}
\date{}
\maketitle
\begin{center}Dedicated to Professor Kiyosi It\^o on his 88th birthday
\end{center}
\begin{abstract}  A Harris flow is a stochastic flow on the real line given by SDE (2.1) below.  We study the {\em noise} generated by Harris flows, 
particularly {\em spectra} of the noise.  Our aim is to understand what  lies beyond the finite order terms in the chaos expansion (the Wiener-It\^o 
expansion) for {\em nonstrong} solutions of SDE (2.1).
\end{abstract}
\setcounter{section}{0}
\section{Definitions and main results}  The notion of noises in continuous time (i.e., the case of time $t\in {\bf R}$) has been introduced by Tsirelson (cf. [T 1], [T 2], [T 5]): 
\begin{defn}  A {\em noise} ${\bf N}=\left[ \{{\cal F}_{s,t}\}_{s\leq t}, \{T_h\}_{h\in {\bf R}}\right]$  is a two parameter family of sub $\sigma$-fields ${\cal F}_{s,t}$, $s\leq t$, of events defined on a probability space $(\Omega, {\cal F}, P)$ which is stationary in time and possesses the following property:
\begin{equation} {\cal F}_{s,u}={\cal F}_{s,t}\otimes {\cal F}_{t,u},\quad s\leq t\leq u,\end{equation}
that is, ${\cal F}_{s,t}$ and ${\cal F}_{t,u}$ are independent and generate ${\cal F}_{s,u}$, for every $s\leq t\leq u$.  By the stationarity in time, we mean the existence of a measurable flow $\{T_h\}$, i.e., a measurable one-parameter group of automorphisms, on $(\Omega,{\cal F}_{-\infty,\infty}:=\bigvee_{s\leq t}{\cal F}_{s,t})$, in which ${\cal F}_{s,t}$ is sent to ${\cal F}_{s+h,t+h}$ by $T_h$.
\end{defn}
In this article, it is always assumed that the probability space is complete and separable and that a sub $\sigma$-field contains all $P$-null sets. 

In the discrete time case (i.e., the case of time $n\in {\bf Z}$), a noise can be defined similarly but it is essentially equivalent to giving an i.i.d. random sequence. In the continuous time case, noises generated by increments of a Wiener process (of finite or countably infinite dimension), a stationary Poisson point process, or an independent pair of them, are typical examples which we call {\em white}, {\em linearizable} or {\em classical noises}.  There are many non-classical noises, however.  Every noise ${\bf N}=\{{\cal F}_{s,t}\}$ contains a unique  maximal (i.e., the largest) classical subnoise which is denoted by ${\bf N}^{lin}=\{{\cal F}_{s,t}^{lin}\}$.

A {\em Harris flow} (as will be defined precisely in Def.1.3 below) is a stochastic flow on the real line ${\bf R}$ determined uniquely by giving a real positive definite function $b(x)$ such that $b(0)=1$, (cf. [H]). Note that $b(x)=b(-x)$. We assume that either $b(x)={\bf 1}_{\{0\}}(x)$ or $b(x)$ is continuous, ${\cal C}^2$ on ${\bf R}\setminus{\{0\}}$ and  {\em strictly positive-definite} in the sense that the matrix $\{b(x_i-x_j)\}$ is strictly positive-definite for any choice of finite different points $\{x_i\}$ in ${\bf R}$. The Harris flow in the discontinuous case of  $b(x)={\bf 1}_{\{0\}}(x)$ is known as the {\em Arratia flow} ([A]).  

Here is a formal definition of stochastic flows on the real line: Let ${\cal T}$ be the set of all non-decreasing right-continuous functions $\varphi: x\in{\bf R}\mapsto \varphi(x)\in {\bf R}$ with the metric defined by  $\rho(\varphi,\psi)=\sum_{n=1}^\infty 2^{-n}\left(\rho_n(\varphi,\psi)\wedge 1\right)$ where
$$\rho_n(\varphi,\psi)=\inf\{ \ \varepsilon>0 \ | \ \varphi(x-\varepsilon)-\varepsilon \leq \psi(x) \leq \varphi(x+\varepsilon)+\varepsilon  \quad \mbox{for all} \ x \in [-n,n] \}.
$$
Then ${\cal T}$ is a Polish space: The composite $(\varphi,\psi)\in {\cal T}\times{\cal T}\mapsto \psi\circ\varphi\in{\cal T}$, defined by $\psi\circ\varphi(x)=\psi(\varphi(x))$, and the evaluation map ${\cal T}\times{\bf R}\ni(\varphi,x)\mapsto \varphi(x)\in{\bf R}$ are all Borel measurable even though they are generally not continuous.
\begin{defn} By a {\em stochastic flow} on ${\bf R}$, we mean a family  ${\bf X}=\{X_{s,t} ; \ s\leq t \ \}$  of ${\cal T}$-valued random variables $X_{s,t}$ having the following properties:
\begin{enumerate}
\item[(1)](Flow property), \ $X_{s,u}=X_{t,u}\circ X_{s,t}$ and $ X_{t,t}={\rm id}$, \  a.s. for every \ $s\leq t \leq u$, 
\item[(2)](Independence property), \ for any sequence $t_0\leq t_1\leq \cdots\leq t_n$,  ${\cal T}$-valued random variables $X_{t_{k-1},t_k}$, $k=1,\cdots,n$, are independent,
\item[(3)](Stationarity), \ for any $h>0$, $X_{s,t}\stackrel{d}{=}X_{s+h,t+h}$,
\item[(4)](Stochastic continuity), \ \ $ X_{0,h}\to{\rm id}$ \ \  in probability as \ \ $ h\downarrow 0.$
\end{enumerate}
\end{defn}
Given a stochastic flow ${\bf X}=\{X_{s,t}\}$, it generates a noise ${\bf N}^X=\left[ \{{\cal F}_{s,t}^X \}, \{T_h\} \right]$ by letting ${\cal F}_{s,t}^X$ to be the $\sigma$-field generated by ${\cal T}$-valued random variables $X_{u,v}$, $s\leq u\leq v\leq t$, and $\{T_h\}$ to be a unique one-parameter family of automorphisms on $(\Omega, {\cal F}_{-\infty,\infty}^X)$  such that  
$(T_h)_*(X_{u,v}(x))=X_{u+h,v+h}(x) $,  $ u\leq v $, $ x\in {\bf R}$. 

Now we give a formal definition of Harris flows.  Generally, for a given filtration ${\bf F}=\{{\cal F}_t\}_{t\geq 0}$, we denote by ${\cal M}_2({\bf F})$ the space of all locally square-integrable ${\bf F}$-martingales $M=(M_t)_{t\geq 0}$ with $M_0=0$ and by ${\cal M}_2^c({\bf F})$ the subspace formed of all continuous elements in ${\cal M}_2({\bf F})$.
\begin{defn}  The {\em Harris flow ${\bf X}=\{X_{s,t}\}$ associated with the correlation function} $b(x)$  is  a stochastic flow on ${\bf R}$ such that, for every $x\in {\bf R}$, if we define the process $M(x)=(M_t(x))_{t\geq 0}$ by setting $M_t(x)=X_{0,t}(x)-x$ and the filtration ${\bf F}^X=\{{\cal F}^X_t\}$ by setting ${\cal F}^X_t={\cal F}_{0,t}^X$, then $M(x)\in {\cal M}_2^c({\bf F}^X)$  and, for every $x, y\in {\bf R}$, we have
\begin{equation} \langle M(x), M(y)\rangle_t =\int_0^t b\left(X_{0,s}(x)-X_{0,s}(y)\right)ds.
\end{equation}
\end{defn}
The law of a Harris flow is uniquely determined under our assumption on functions $b(x)$.  The existence of Harris flows has been established in [H] (cf. also [LR 1]).  A Harris flow is equivalently given by a stochastic differential equation (SDE) (2.1) in Section 2.

Let ${\bf X}=\{X_{s,t}\}$ be a Harris flow associated with the function $b(x)$ and ${\bf N}^X$ be the noise generated by it. {\em Suppose that $b(x)$ is continuous}. Then we can construct a centered Gaussian system  ${\bf W}= \{W(t,x); t\in{\bf R}, x\in{\bf R}\}$ contained in $L_2({\cal F}_{-\infty,\infty}^X)$ such that $(T_h)_*[W(t,x)-W(s,x)]=W(t+h,x)-W(s+h,x)$, $s\leq t$, $x\in {\bf R}$  and, if we set $w_t(x)=W(t,x)-W(0,x)$, then $w(x)=(w_t(x))_{t\geq 0}\in {\cal M}_2^c({\bf F}^X)$  and, for every $x, y\in {\bf R}$, we have 
$\langle w(x), w(y)\rangle_t = t b(x-y)$.  Indeed, $W(t,x)-W(s,x)$ is the $L_2$-limit of $M_\Delta^x(s,t)$ as $|\Delta|\to 0$.  Here, for a sequence of times $\Delta: s=t_0<t_1<\cdots<t_{n-1}<t_n=t$ and $x\in {\bf R}$, $M_\Delta^x(s,t)=\sum_{k=1}^n (X_{t_{k-1},t_k}(x)-x)$ and $|\Delta|=\max_k|t_k-t_{k-1}|$. ${\bf W}$ defines a Gaussian white noise  ${\bf N}^W=\left[\{{\cal F}^W_{s,t}\}, \{T_h\}\right]$ where ${\cal F}^W_{s,t}=\sigma[W(v,x)-W(u,x); s\leq u\leq v\leq t, x\in{\bf R}]$. It is obvious that ${\bf N}^W$ is a subnoise of  ${\bf N}^X$.
 \begin{th} Suppose that the function $b(x)$ is continuous. Then, it holds that
 $[{\bf N}^X]^{lin}={\bf N}^W $.  Furthermore,  ${\bf N}^X={\bf N}^W $  holds,
that is, the noise ${\bf N}^X$ generated by the Harris flow ${\bf X}$ is classical, if and only if
\begin{equation} \int_{0+}^1(1-b(x))^{-1}dx =\infty.\end{equation}
\end{th}
Hence, the noise ${\bf N}^X$ is nonclassical if and only if
\begin{equation} \int_{0+}^1(1-b(x))^{-1}dx <\infty.\end{equation}

In the case of the Arratia flow, it generates a nonclassical noise: Tsirelson [T 3] (cf. also [LR 2]) showed that this noise is {\em black} in the sense that $({\cal F}^X_{s,t})^{lin}=\{\emptyset, \Omega\}$ for every $s\leq t$.

Tsirelson ([T 2], [T 5]) introduced the notion of {\em spectral measures} for noises which is an invariant under the isomorphism of noises and which can measure the degree of non-linearity (or sensitivity in the discrete-time approximation) of noises. Let ${\cal C}$ be the space formed of all compact sets in ${\bf R}$ endowed with the Hausdorff distance and ${\cal C}^f$ be its subclass formed of all finite sets: ${\cal C}^f=\{ \ S\in {\cal C} \ | \  |S| <\infty \ \}$. Here, $|S|$ denotes the number of elements in $S$.
\begin{defn}  Let ${\bf N}=\left[\{{\cal F}_{s,t}\}, \{T_h\} \right]$ be a noise. To every $\Phi\in L_2({\cal F}_{-\infty,\infty})$, there corresponds a unique finite Borel measure $\mu_\Phi$ on ${\cal C}$ such that
\begin{equation} \mu_\Phi( \ \{S\in {\cal C} \ | \ S \subset J\} \ )=E\left[E(\Phi|{\cal F}(J))^2\right] \end{equation}
for every elementary set $J\subset {\bf R}$.  Here, by an elementary set $J$, we mean a finite union $J=\bigcup_k [t_{2k},t_{2k+1}]$ of non-overlapping intervals and we set ${\cal F}(J)=\bigvee_k {\cal F}_{t_{2k},t_{2k+1}}$.  $\mu_\Phi$ is called the {\em spectral measure} of the noise ${\bf N}$ associated with $\Phi\in L_2({\cal F}_{-\infty,\infty})$.
\end{defn}
When  $\Phi\in L_2({\cal F}_{s,t})$,  we have $\mu_\Phi({\cal C}\setminus{\cal C}_{[s,t]})=0$ where ${\cal C}_{[s,t]}=\{S\in {\cal C} \ | \ S \subset [s,t]\}$, so that $\mu_\Phi$ is a measure on ${\cal C}_{[s,t]}$.  The following is an important characterization of classical noises due to Tsirelson: {\em a noise is classical if and only if $\mu_\Phi({\cal C}\setminus{\cal C}^f)=0$ for every $\Phi\in L_2({\cal F}_{-\infty,\infty})$.}

Set $L_2^{us}({\cal F}_{s,t})=\{ \ \Phi\in L_2({\cal F}_{s,t}) \ | \ ||\Phi||_2=1 \ \}$; the unit sphere in $L_2({\cal F}_{s,t})$. If $\Phi\in L_2^{us}({\cal F}_{-\infty,\infty})$, then $\mu_\Phi$ is a Borel probability on ${\cal C}$ so that we can speak of a ${\cal C}$-value random variable with the distribution $\mu_\Phi$.  We denote it by $S_\Phi$ and call it the {\em spectral set} of the noise associated with $\Phi$.

We wish to describe the spectral set $S_\Phi$ for the noise ${\bf N}^X$ generated by a Harris flow ${\bf X}$ when $\Phi=X_{0,1}(0)\in L_2^{us}({\cal F}^X_{0,1})$.  The random set $S_\Phi$ in this case is denoted by $S_X$.  We would also obtain some information on $S_\Phi$ for general $\Phi$. We consider naturally the case when the noise is nonclassical so that we assume (1.4). Furthermore, we assume that 
\begin{equation} b(x)\quad \mbox{\em is non-increasing in} \ (0,\infty) \ \mbox {\em and satisfies} \ \  \lim_{x\to \infty}b(x)=0.\end{equation}
Functions $b(x)=\exp(-c|x|^\a)$ for $c>0$ and $0<\a<1$ are typical examples. Also, $b(x)={\bf 1}_{\{0\}}(x)$ (the case of the Arratia flow) is another typical example.

For $S\in {\cal C}$, let $S^{acc}$ be the the set of all accumulation points of $S$, so that $S^{acc}\ne \emptyset$ if and only if $S\notin {\cal C}^f$.
\begin{th} 
 Let ${\bf X}$ be the Harris flow associated with the function $b(x)$ which satisfies (1.4) and (1.6) and let $S_X$ be the spectral set $S_\Phi$ of the noise ${\bf N}^X$ for $\Phi=X_{0,1}(0)$. Then the random set $S_X^{acc}$ has the same law as the random set ${\widetilde S}$ in $[0,1]$ defined by
\begin{equation}{\widetilde S}=\{ \ t \ | \ 0\leq t\leq \tau, \ {\widehat \xi}^+(\tau-t)=0 \ \}\end{equation}
 where ${\widehat \xi}^+=\{{\widehat \xi}^+(t)\}_{t\geq0}$ is the reflecting diffusion process on $[0,\infty)$ with the generator 
\begin{equation}{\widehat L}=\frac{d}{dx} (1-b(x))\frac{d}{dx}\end{equation}
and the initial distribution $\mu(dx):=-db(x)$.  Here, $\tau$ is a $[0,1]$-valued and uniformly distributed random variable independent of ${\widehat \xi}^+$.
\end{th}
In particular, we have
$$P(S_X^{acc}\ne\emptyset)=P(|S_X|=\infty)=P({\widetilde S}\ne \emptyset)= P\left\{ \exists t\in [0,\tau]; {\widehat \xi}^+(t)=0 \right\}$$
and this probability is also equal to $ E\left[\int^1_0(1-b(\xi^+(t)))dt\right]$ where $\xi^+=\{\xi^+(t)\}_{t\geq 0}$ is the reflecting diffusion process on $[0,\infty)$ with the generator 
\begin{equation} L=(1-b(x))\frac{d^2}{dx^2}\end{equation}
which starts at $0$.  Still another expression of this probability is given by the expectation 
$\frac12 E\left[A^{-1}(1)\right]$,  where $A(t)$ is an additive functional of the one-dimensional Wiener process $\b(t)$ with $\b(0)=0$, defined by
\begin{equation} A(t)=\frac12 \int^t_0(1-b(\b(s)))^{-1} ds, \end{equation}
and $t\to A^{-1}(t)$ is the inverse function of $t\to A(t)$.

In the case of the Arratia flow, $S_X^{acc}=S_X$ and it is a perfect set, a.s.. It is described as a zero points set of a (double speed) reflecting Brownian motion starting at $0$ as in the theorem. This recovers a result of Tsirelson ([T 4]) who obtained it by an approximation by coalescing random walks.

In the following, we consider the class of Harris flows associated with the correlation functions $b(x)$ which satisfy (1.4), (1.6) and, for some $0\leq \a<1$,
\begin{equation}   1-b(x)\asymp |x|^{\a}\quad\mbox{as} \ \ x\to 0.\end{equation}
Again, functions  $b(x)=\exp(-c|x|^\a)$ for $c>0$ and $0<\a<1$ are typical examples.  Note also that the function  $b(x)={\bf 1}_{\{0\}}(x)$ (the case of the Arratia flow) is a typical example of the case when $\a=0$. From Theorem 1.2, we can obtain the following: Denoting by $\dim(S)$ the Hausdorff dimension of a subset $S$ in ${\bf R}$,
\begin{cor}
$\dim(S_X^{acc})=\frac{1-\a}{2-\a} \quad {\rm a.s.}$, under the condition that it is not empty.
\end{cor}
\begin{th} Let $\gamma=\inf\{ \ \b \ | \ \dim(S_\Phi)\leq \b, \ {\rm a. s.} \ \mbox{for any} \ \Phi\in L_2^{us}({\cal F}^X_{-\infty,\infty}) \ \}$.
Then
$$ \gamma = \frac{1-\a}{2-\a}.$$
\end{th}

The proof of these theorems will be given in the subsequent sections by appealing to two main tools: {\em joinings of Harris flows} and  certain {\em duality relations} between the reflecting (absorbing) $L$-diffusion and the absorbing (resp. reflecting) ${\widehat L}$-diffusion.  
\setcounter{section}{1}
\setcounter{equation}{0}
\section{The joining of Harris flows: The proof of Th. 1.1.}
Suppose that the correlation function $b(x)$ of a Harris flow ${\bf X}$ is continuous. Let $H (\subset {\bf C}_b({\bf R}\to {\bf R}) )$ be the (real) reproducing kernel Hilbert space associated with $b(x)$ so that, defining $f_x\in H$ by $f_x(y)=b(y-x)$,  linear combinations $\sum c_i f_{x_i}$ are dense in $H$ and $(f_x, f_y)_H=b(x-y)$. The Gaussian system ${\bf W}$ introduced in Section 1 can be given equivalently by a Gaussian system $\{W(t,f); t\in{\bf R}, f\in H\}$ contained in $L^2({\cal F}^X_{-\infty,\infty})$  such that $(T_h)_*[W(t,f)-W(s,f)]=W(t+h,f)-W(s+h,f)$, $s\leq t$, $f\in H$  and, if we set $w_t(f)=W(t,f)-W(0,f)$, then $w(f)=(w_t(f))_{t\geq 0}\in {\cal M}_2^c({\bf F}^X)$  and, for every $f, g\in H$, we have 
$\langle w(f), w(g)\rangle_t = t(f,g)_H$. Indeed, we set $W(t,f)=\sum_i c_i W(t,x_i)$ when $f=\sum c_i f_{x_i}$ and extend this to general $f\in H$ by routine arguments.

We define an It\^o-type stochastic integral  $\int_0^t \psi_s \cdot W(ds,\varphi_s)$  for ${\bf F}^X$-predictable processes $\varphi$ and $\psi$ satisfying that $\int^t_0|\psi_s|^2ds<\infty$, a.s., by
$$\int_0^t \psi_s \cdot W(ds,\varphi_s)=\sum_k \int^t_0\psi_s\cdot e_k(\varphi_s)db_k(s),$$
where $\{e_k\}$ is an orthonormal basis (ONB) in $H$ and $b_k(t)=W(t,e_k)$, so that $\{b_k(t)\}$ is an independent family of one-dimensional Wiener processes. As is easily seen, the definition is independent of a particular choice of ONB.
Note that $\sum_k e_k(\varphi_s)e_k(\varphi'_s)=b(\varphi_s-\varphi'_s)$, so that, in particular, $\sum_k |e_k(\varphi_s)|^2\equiv 1$. Now, (1.2) is equivalently given in the form of SDE for $X_t:=X_{0,t}(x)$:
\begin{equation}
X_t=x+\int_0^t W(ds,X_s)=x+\sum_k \int^t_0 e_k(X_s)db_k(s).
\end{equation}
Since  $\sum_k |e_k(x)-e_k(y)|^2=2(1-b(x-y))$,  the condition (1.3) implies the pathwise uniqueness of solutions for SDE (2.1) (cf. [IW], p.182).  Hence, if the function  $b$  satisfies the condition (1.3), then $X_t$  is a unique strong solution to SDE (2.1) so that $X_{0,t}(x)$ is ${\cal F}^W_{0,t}$-measurable for every $x$.  By the stationarity, we see that $X_{s,t}(x)$ is ${\cal F}^W_{s,t}$-measurable for every $x$ and $s\leq t$. Therefore, $ {\bf N}^X={\bf N}^W$  holds.
Thus, the {\em if part} of Th. 1.1 is proved.

To prove the {\em only if part}, we first remark the following martingale representation theorem for Harris flows.
\begin{prop}  Suppose the correlation function $b(x)$ of the Harris flow is continuous. Then,
$M\in {\cal M}_2({\bf F}^X)$   if and only if there exists a sequence  $\varphi_k=(\varphi_k(t))$, $k=1,2,\ldots,$ of ${\bf F}^X$-predictable processes satisfying that  $\sum_k\int^t_0\varphi_k^2(s)ds<\infty$, a.s.,for each  $t>0,$
and
$$ M(t)=\sum_k\int^t_0\varphi_k(s)db_k(s).$$
In particular, it holds that ${\cal M}_2({\bf F}^X)={\cal M}_2^c({\bf F}^X).$
\end{prop}
{\em Proof.} \ \ 
Given distinct $x_1, x_2, \ldots x_n \in {\bf R}$, any ${\bf R}^n$-valued process $(X^1_t, X^2_t,\ldots X^n_t)$  of which each component $X^k_t$ solves the SDE (2.1) starting from $x_k$ and these components satisfy the coalescing property, has the same law as the $n$-point motion of the Harris flow  \\ $(X_{0,t}(x_1), X_{0,t}(x_2) \ldots , X_{0,t}(x_n))$. From this uniqueness in law, it follows by the usual methods that any $M \in {\cal M}_2({\bf F}^X)$ that is measurable with respect to this $n$-point motion is continuous and has the desired representation as a stochastic integral. The result can then be extended to an arbitrary $M \in {\cal M}_2({\bf F}^X)$ using the fact that the set of representable martingales is closed in this space.
 \hfill \qed \medskip

>From this proposition, we can easily deduce that $[{\bf N}^X]^{lin}={\bf N}^W $, see also Lemma 6a5 of [T 5].  Indeed, if ${\bf N}^W$ is smaller than $[{\bf 
N}^X]^{lin}$, then there should exist some martingale in ${\cal M}_2({\bf F}^X)$ which cannot be given by a sum of stochastic integrals by $b_k$.  Hence, in 
order to prove the {\em only if part}, it is sufficient to show that {\em (1.4) implies that ${\bf N}^W$ is strictly smaller than ${\bf N}^X$}. For 
this, we introduce the following notion.

\begin{defn} By a {\em joining of a Harris flow}, we mean a pair $( {\bf X}=\{X_{s,t}\}, {\bf X}'=\{X'_{s,t}\} )$ of copies of the Harris flow defined on a same probability space such that the joint process $\Xi=\{\Xi_{s,t}=(X_{s,t},X'_{s,t});s\leq t\}$ has the independence property (2) in Def.1.2.  Given $0\leq \rho\leq 1$, it is called a {\em $\rho$-joining} if it satisfies further the following: ${\bf X}$ and ${\bf X}'$ are stationarily correlated in the sense that the joint process $\Xi$ has the stationarity 
property (3)  of  Def.1.2 and, if filtrations ${\bf F}^X=\{{\cal F}^X_t\}$, ${\bf F}^{X'}=\{{\cal F}^{X'}_t\}$ and martingales $M(x)=(M_t(x))$, $M'(x)=(M'_t(x))$ are defined similarly as in Def.1.3 for ${\bf X}$ and ${\bf X}'$, respectively, then  ${\bf F}^X$ and ${\bf F}^{X'}$ are jointly immersed, i.e., ${\cal M}_2({\bf F}^X)\cup{\cal M}_2({\bf F}^{X'}) \subset {\cal M}_2({\bf F}^X\bigvee{\bf F}^{X'})$, and, for every $x, y\in {\bf R}$,
\begin{equation} \langle M(x), M'(y)\rangle_t =\int_0^t \rho\cdot b\left(X_{0,s}(x)-X'_{0,s}(y)\right)ds,
\end{equation}
$b(x)$ being the correlation function of the Harris flow.
\end{defn}
It is obvious that, for a $\rho$-joining, the corresponding Gaussian noises ${\bf W}$ and ${\bf W}'$ are jointly Gaussian and $\rho$-correlated.
\begin{lm}
For  $0\leq\rho<1$, a $\rho$-joining exists and is unique in law. If, in particular, $\rho=0$, then it is a pair of independent copies.
\end{lm}
This lemma can be deduced from the fact that the following differential operator $\Lambda$ with variables $x=(x_1,\cdots,x_n)\in {\bf R}^n$ and $x'=(x'_1,\cdots,x'_m)\in {\bf R}^m$ is non degenerate at all such points $(x,x')\in {\bf R}^n\times{\bf R}^m$ as all coordinates in $x$ are different and also all coordinates in $x'$ are different:
\begin{eqnarray*}\Lambda&=&\frac12\sum_{i=1}^n \sum_{j=1}^n b(x_i-x_j)\frac{\partial^2}{\partial x_i \partial x_j}+\frac12\sum_{k=1}^m \sum_{l=1}^m b(x'_k-x'_l)\frac{\partial^2}{\partial x'_k \partial x'_l}\\
&+&\rho\sum_{i=1}^n\sum_{k=1}^m b(x_i-x'_k)\frac{\partial^2}{\partial x_i \partial x'_k}.
\end{eqnarray*}
Note that, for a $\rho$-joining $({\bf X}, {\bf X}')$, the process
$$[0,\infty)\ni t\mapsto (X_{0,t}(x_1),\cdots, X_{0,t}(x_n), X'_{0,t}(x'_1),\cdots, X'_{0,t}(x'_m))$$
is a solution to the $\Lambda$-martingale problem. \medskip

We now assume (1.4) and prove that ${\bf N}^W$ is strictly smaller than ${\bf N}^X$. Take  $\rho$-joinings $({\bf X}^{(\rho)}, {\bf X}'^{(\rho)})$  for $\rho\in[0,1)$.  By (2.2), the process $\xi^{(\rho)}(t)=X_{0,t}^{(\rho)}(0)-X_{0,t}'^{(\rho)}(0)$ is a Feller diffusion on ${\bf R}$ with the canonical scale $s(x)=x$ and the speed measure $m(dx)=(1-\rho\cdot b(x))^{-1}dx$ which starts from the origin at time $0$, (cf. [IM] for a general theory of Feller diffusions).  As $\rho\nearrow 1$, the processes $\xi^{(\rho)}(t)$  converge to the Feller diffusion $\xi(t)$  with the canonical scale $s(x)=x$ and the speed measure $m(dx)=(1- b(x))^{-1}dx$ which starts from the origin $0$ at time $0$.  As is well-known, $\xi(t)=\b(A^{-1}(t))$  for a one-dimensional Wiener process $\b(t)$ and $A(t)$ is defined by (1.10). Then we have 
$$\lim_{\rho\nearrow 1}E\left[|\xi^{(\rho)}(t)|^2\right]=E\left[|\xi(t)|^2\right]=\frac12 E[A^{-1}(t)]>0$$
for $t>0$.  Suppose ${\bf N}^X\subset{\bf N}^W$ be true.  Then $X_{0,t}^{(\rho)}(0):=\Phi\in L_2({\cal F}_{0,t}^W)$ and $E\left[X_{0,t}'^{(\rho)}(0)|{\bf W}\right]=P_{-\log \rho}\Phi$ where $(P_s)_{s\geq 0}$ is the {\em Ornstein-Uhlenbeck semigroup} acing on $L_2({\cal F}^W_{0,t})$. \\
 Hence $E\left[|\xi^{(\rho)}(t)|^2\right]=2\left(||\Phi||_2^2-(\Phi,P_{-\log \rho}\Phi)_2||\right)$. By the $L^2$-continuity of the Ornstein-Uhlenbeck semigroup, we have
$$\lim_{\rho\nearrow 1}E\left[|\xi^{(\rho)}(t)|^2\right]=\lim_{\rho\nearrow 1} 2\left(||\Phi||_2^2-(\Phi,P_{-\log \rho}\Phi)_2||\right)=0.$$
Thus we have a contadiction and hence we cannot have ${\bf N}^X\subset{\bf N}^W$. This proves the {\em only if part} of Th.1.1 so that its proof now is completed.\medskip

In the following, we assume that (1.4) holds so that the noise generated by the Harris flow is nonclassical.  In this case, $1$-joinings are not unique. We specify two of them as the {\em $1^+$-joining} and the {\em $1^-$-joining}.
\begin{defn} The $1^+$-joining $({\bf X}, {\bf X}')$  is the identity joining: i.e., ${\bf X}={\bf X}'$. The  $1^-$-joining is the limit in law
of the  $\rho$-joinings $({\bf X}^{(\rho)}, {\bf X}'^{(\rho)})$ as $\rho\nearrow 1$ .  It is  such that $[0,\infty)\ni t\mapsto X_{0,t}(x)-X'_{0,t}(y)$, for 
fixed $x, y\in {\bf R}$, is the Feller diffusion on ${\bf R}$ with the canonical scale $s(x)=x$ and the speed measure $m(dx)=(1-b(x))^{-1}dx$ which starts at 
$x-y$ at time $0$. \end{defn}
For $\rho\in [0,1)$, let  $({\bf X}, {\bf X}^\prime)$ be a $\rho$-joining with  corresponding  $\rho$-correlated Gaussian processes   ${\bf W}$ and ${\bf W}^\prime$.  It is easy to see that the joint law $\Pi(d{\cal X}d{\cal X}^\prime d{\cal W}d{\cal W}^\prime)$ of $({\bf X}, {\bf X}^\prime, {\bf W},{\bf W}^\prime)$ is given by
$$ P({\bf X}\in d{\cal X}|{\bf W}={\cal W})P({\bf X}^\prime\in d{\cal X}^\prime|{\bf W}^\prime={\cal W}^\prime)P({\bf W}\in d{\cal W},{\bf W}^\prime\in d{\cal W}^\prime).$$
From this, we deduce that
\begin{eqnarray*}
 & E[ \Phi\cdot \pi_*(\Psi)] & = E\bigl[ E[ \Phi | {\bf W} ]\cdot E[ \pi_*(\Psi)|  {\bf W}^\prime] \bigr] \\
& &=E\bigl[ E[ \Phi | {\bf W} ] \cdot E\bigl(E[ \pi_*(\Psi)|  {\bf W}^\prime] |{\bf W}\bigr) \bigr]\\
& &=E\bigl[ E[ \Phi | {\bf W} ] \cdot E\bigl[\pi_*(E(\Psi|{\bf W}) |{\bf W}\bigr]\bigr]=E\bigl[E[ \Phi | {\bf W} ]\cdot P_{-\log \rho}(E(\Psi|{\bf W}))\bigr]
\end{eqnarray*}
whenever $\Phi, \Psi \in L_2({\cal F}^X_{-\infty, \infty})$. Here, $\pi_*$ is the unique isomorphism \\
 $\pi_*: L_2({\cal F}^X_{-\infty, \infty})\to  L_2({\cal F}^{X'}_{-\infty,\infty})$  such that $\pi_*(X_{s,t}(x))=X'_{s,t}(x)$ for every $s,t$ and $x$, and $(P_s)$ is the Ornstein-Uhlenbeck semigroup acting on $L_2({\cal F}^W_{-\infty, \infty})$.  By the $L^2$-continuity of the Ornstein-Uhlenbeck semigroup, the above expectation converges to $E\bigl[E[ \Phi | {\bf W} ]\cdot E[ \Psi|  {\bf W}]\bigr]$  as $\rho\nearrow 1$.  This proves existence of the $1^-$-joining as the limit of $\rho$-joinings.  Moreover  for a $1^-$-joining $({\bf X}, {\bf X}^\prime)$ the corresponding Gaussian systems  ${\bf W}$ and ${\bf W}^\prime$ are equal and  ${\bf X}$ and ${\bf X}^\prime$ are conditionally independent given this common Gaussian process.  
\begin{rem}
For the Arratia flow, its $\rho$-joining for $\rho\in[0,1)$ is independent of $\rho$ and coincides with $0$-joining, that is, a pair of independent copies of the Arratia flow. Hence, its $1^-$-joining is also a pair of independent copies of the Arratia flow.
\end{rem}

   Let $F=\bigcup_{k=1}^n[t_{2k-2},t_{2k-1}]$ be an elementary set in ${\bf R}$ defined for a sequence $t_0<t_1<\cdots<t_{2n-2}<t_{2n-1}$ of times. We would introduce the notion of {\em $(\rho, F)$-joining} $({\bf X}, {\bf X}')$ of the Harris flow when $\rho\in[0,1)$, which is roughly the $\rho$-joining on $F$ and the identity joining outside $F$. To be more precise, set $t_{-1}=-\infty$ and $t_{2n}=\infty$ by convention. Take a $\rho$-joining $({\bf Y}, {\bf Y}')$ and a $1^+$-joining $({\bf Z}, {\bf Z}')$ which are mutually independent.  Define ${\bf X}=[\{X_{s,t}\}_{s\leq t}]$ as follows:  First, set $X_{s,t}=Y_{s,t}$  if $t_{2k-2}\leq s\leq t\leq t_{2k-1}$, $k=1,\cdots,n$ and $X_{s,t}=Z_{s,t}$ if $t_{2k-1}\leq s\leq t\leq t_{2k}$, $k=0,\cdots,n$.  Then, define $X_{s,t}$ for general $ s\leq t$, by
$$ X_{s,t}=X_{t_l,t}\circ X_{t_{l-1},t_l}\circ\cdots\circ X_{t_k,t_{k+1}}\circ X_{s,t_k}$$
when $ t_{k-1}<s\leq t_k\leq t_l\leq t<t_{l+1}$, $0\leq k\leq l\leq 2n-1$. Define ${\bf X}'=[\{X'_{s,t}\}_{s\leq t}]$ similarly from ${\bf Y}'$ and ${\bf Z}'$. Then $({\bf X}, {\bf X}')$ defines a joining of the Harris flow in which, however, ${\bf X}$ and ${\bf X}'$ are not stationarily correlated.
\begin{defn}
The pair $({\bf X}, {\bf X}')$ defined above is called the {\em $(\rho, F)$-joining} of the Harris flow.
\end{defn}
Next, take mutually independent $1^-$-joining $({\bf Y}, {\bf Y}')$ and $1^+$-joining $({\bf Z}, {\bf Z}')$ and construct the pair $({\bf X}, {\bf X}')$ in the same way.
\begin{defn}
The pair $({\bf X}, {\bf X}')$ defined above is called the {\em $(1^-, F)$-joining} of the Harris flow
\end{defn}

We turn now to the notion of the  spectral measure $\mu_\Phi$ associated with  some $ \Phi\in L_2({\cal F}^X_{-\infty,\infty})$ as defined in Def.1.4.  This notion is intimately 
related to chaos expansions.  The spectral measure of a random variable $ \Phi\in L_2({\cal F}^W_{-\infty,\infty})$, measurable with respect to ${\bf W}$, can be expressed  by expanding $\Phi$ as a sum of multiple Wiener-It\^o integrals with respect to the Brownian motions $b_k$. To be more precise, $\Phi=\sum_{m=0}^\infty I_m$  where $I_0$ is a constant and $I_m$, for $m=1,2.\cdots,$ is given by
an iterated It\^o stochastic integral
$$ I_m= \sum_{(k_1,\cdots,k_m)}\int\cdots\int_{-\infty<t_m<\cdots<t_1<\infty} f_\Phi^{(k_1,\cdots,k_m)}(t_1,\cdots,t_m)db_{k_m}(t_m)\cdots db_{k_1}(t_1).$$
 $\mu_\Phi$ is supported on ${\cal C}^f=\{S\in {\cal C}: |S|<\infty\}$ and 
\begin{eqnarray*}&  &\mu_\Phi({\cal C}^f)=E(\Phi^2)=\sum_{m=0}^\infty E(|I_m|^2)\\
&=&\sum_{m=0}^\infty \sum_{(k_1,\cdots,k_m)}\int\cdots\int_{-\infty<t_m<\cdots<t_1<\infty} |f_\Phi^{(k_1,\cdots,k_m)}(t_1,\cdots,t_m)|^2 dt_m\cdots dt_1<\infty.\end{eqnarray*}
The restriction of $\mu_\Phi$ to $\{S\in 
{\cal C}:|S|=m\}$ is determined (denoting $S=\{t_m,\cdots,t_1\}$, $-\infty<t_m<\cdots<t_1<\infty$) by
$$\mu_\Phi(dS;|S|=m)=|f_\Phi^{(k_1,\cdots,k_m)}(t_1,\cdots,t_m)|^2 dt_m\cdots dt_1.$$
In particular, $\mu_\Phi(|S|=m)=E(|I_m|^2)$.

For a general $ \Phi\in L_2({\cal F}^X_{-\infty,\infty})$, the chaos expansion of \ $E[ \Phi| {\bf W}]$ \ given by $E[ \Phi| {\bf W}]=\sum_{m=0}^\infty I_m$,  yields in the same fashion  the restriction of $\mu_\Phi$ to ${\cal C}^f$ and in particular 
$$E\left[ E[ \Phi| {\bf W}]^2\right] = \mu_\Phi\bigl({\cal C}^f \bigr). $$

If $({\bf X},{\bf X}^\prime)$ is a $\rho$-joining for $\rho\in[0,1)$ and $\Phi^\prime=\pi_*(\Phi)$ as above, we have
$$ E\bigl(E[ \Phi^\prime| {\bf W}^\prime]|{\bf W}\bigr)=P_{-log \rho}\bigl(E[\Phi|{\bf W}]\bigr)=\sum_{m=0}^\infty\rho^m I_m.$$
As was remarked above, the relation $E(\Phi\Phi^\prime)=E\bigl(E[ \Phi| {\bf W}]E[ \Phi^\prime| {\bf W}^\prime]\bigr)$ holds. Hence, 
\begin{equation}
E\bigl(\Phi\Phi^\prime\bigr)=\sum_{m=0}^\infty \rho^m E(|I_m|^2)=\int_{{\cal C}}\rho^{|S|}\mu_\Phi(dS).
\end{equation}
In the same way, we deduce for a $1^-$-joining $({\bf X},{\bf X}^\prime)$,
\begin{equation}
E\bigl(\Phi \Phi^\prime\bigr) = \mu_\Phi({\cal C}^f).
\end{equation}
\begin{ex} \ \ Consider the case $\Phi=g(X_{0,1}(x))$ for a bounded continuous function $g$ on ${\bf R}$. Note that $E(\Phi^2)=\int_{{\bf R}}p(t,x-y)g(y)^2dy$ where 
$$p(t,x)=\frac{1}{\sqrt{2\pi t}}\exp\{-\frac{x^2}{2t}\}, \quad t>0, \ x\in {\bf R}.$$  
The chaos expansion of $E[\Phi|{\bf W}]$ was obtained explicitly by Veretennikov and Krylov (cf. [VK]): By setting 
$$T_tf(x)=\int_{{\bf R}}p(t,x-y)f(y)dy \quad \mbox{and}\quad Q_t^kf(x)=e_k(x)\frac{\partial}{\partial x}T_tf(x),$$
we have
$$g(X_{0,1}(x))=\sum_{m=0}^n I_m + R_n,\quad I_0=T_1g(x)=E[\Phi],$$
where $I_m$,  $m=1,\ldots,n$, and $R_n$ are given by the following iterated It\^o stochastic integrals:
\begin{eqnarray*} I_m&=& \sum_{(k_1,k_2,\cdots,k_m)}\int\cdots\int_{0<t_m<t_{m-1}<\cdots<t_2<t_1<1} \left[T_{t_m}Q^{k_m}_{t_{m-1}-t_m}\cdots\right. \\
&\cdots&\left.Q^{k_2}_{t_1-t_2}Q^{k_1}_{1-t_1}g(x)\right] db_{k_m}(t_m)db_{k_{m-1}}(t_{m-1})\cdots db_{k_2}(t_2)db_{k_1}(t_1),
\end{eqnarray*}

\begin{eqnarray*} R_n&=& \sum_{(k_1,k_2,\cdots,k_n,k_{n+1})}\int\cdots\int_{0<t_{n+1}<t_n<\cdots<t_2<t_1<1} \left[Q^{k_{n+1}}_{t_n-t_{n+1}}Q^{k_n}_{t_{n-1}-t_n}\cdots\right.\\
&\cdots& \left.Q^{k_2}_{t_1-t_2}Q^{k_1}_{1-t_1}g(X_{0,t_{n+1}}(x))\right] db_{k_{n+1}}(t_{n+1})db_{k_n}(t_n)\cdots db_{k_2}(t_2)db_{k_1}(t_1).
\end{eqnarray*}
>From this, we obtain that
$$E[\Phi|{\bf W}]=\sum_{m=0}^\infty I_m.$$
\end{ex}

The following is a key lemma for the proof of Theorem 1.2  which records various generalizations of the identities (2.3) and (2.4). As above, we denote by $S_X$  the spectral set $S_\Phi$ when $\Phi=X_{0,1}(0)$ which is a ${\cal C}_{[0,1]}$-valued random variable.

\begin{lm}
(i) \  If $({\bf X}, {\bf X}')$ is a $(\rho, F)$-joining of the Harris flow for  $\rho\in[0,1)$, then,
\begin{equation}
E\left[\rho^{|S_X\cap F|}\right]=E\left[X_{0,1}(0)X'_{0,1}(0)\right],
\end{equation}
equivalently,
\begin{equation}
 E\left[1-\rho^{|S_X\cap F|}\right]=\frac12 E\left[|X_{0,1}(0)-X'_{0,1}(0)|^2\right].
\end{equation}
(ii)  \ \ If $({\bf X}, {\bf X}')$ is a $(1^-, F)$-joining of the Harris flow, then,
\begin{equation} P\left( \ |S_X\cap F|<\infty \ \right)=E\left[X_{0,1}(0)X'_{0,1}(0)\right],
\end{equation}
equivalently,
\begin{equation}
P\left( \ |S_X\cap F|=\infty \ \right)=\frac12 E\left[|X_{0,1}(0)-X'_{0,1}(0)|^2\right].
\end{equation}
(iii)  \ \ More generally, let $({\bf X}, {\bf X}')$ be a $(\rho, F)$-joining for $0\leq\rho<1$ (a $(1^-, F)$-joining) and $\Phi\in L_2^{us}({\cal F}^X_{-\infty,\infty})$.  There is a unique isomorphism $\pi_*: L_0({\cal F}^X_{-\infty, \infty})\to  L_0({\cal F}^{X'}_{-\infty,\infty})$  such that $\pi_*(X_{s,t}(x))=X'_{s,t}(x)$ for every $s,t$ and $x$.  Set $\Phi'=\pi_*(\Phi)$.  Then we have
\begin{equation} E\left[\rho^{|S_\Phi\cap F|}\right] \ \bigl( \ {\rm resp.} \  P\left( \ |S_\Phi\cap F|<\infty \ \right) \ \bigr) =E\left[\Phi \Phi' \right],
\end{equation}
equivalently,
\begin{equation}  E\left[1-\rho^{|S_\Phi\cap F|}\right] \ \bigl( \ {\rm resp.} \ P\left( \ |S_\Phi\cap F|=\infty \ \right) \ \bigr)=\frac12 E\left[|\Phi-\Phi'|^2 \right]\end{equation}
\end{lm}
{\em Proof}.  \ \ In the case when $\Phi\in L_2^{us}({\cal F}^X_{s,t})$ and $F=[s,t]$, (2.9) is nothing but (2.3) and (2.4).  From this, we can deduce (2.9) in the general case of an elementary set $F=\bigcup_{k=1}^n [t_{2k-2},t_{2k-1}]$, $t_{-1}=-\infty<t_0<\cdots<t_{2n-1}<t_{2n}=\infty$, by considering the following  $L^2$-space factorization:
$$ L_2({\cal F}_{-\infty,\infty}^X)=\bigotimes_{k=0}^{2n} L_2({\cal F}_{t_{k-1},t_k}^X).$$ 
We omit the details. \hfill \qed
\setcounter{section}{2}
\setcounter{equation}{0}
\section{Duality relations for $L$- and ${\widehat L}$-diffusions in the time reversal: The proof of Th. 1.2.}  
Let $\{\xi^+(t), P_x\}$ and $\{{\widehat \xi}^+(t), {\widehat P}_x\}$ be the reflecting $L$- and ${\widehat L}$-diffusion processes on $[0,\infty)$ introduced in Section 1.  The associated Markovian semigroups of operators acting on the space ${\bf B}([0,\infty))$ of real bounded Borel functions are defined by
\begin{equation}
T^+_t f(x)=E_x[f(\xi^+(t))] \quad \mbox{and}\quad {\widehat T}^+_t f(x)={\widehat E}_x[f({\widehat \xi}^+(t))].
\end{equation}
Define also the semigroups for {\em absorbing processes} by
\begin{equation}
T^-_t f(x)=E_x[f(\xi^+(t\wedge\sigma_0))] \quad \mbox{and}\quad {\widehat T}^-_t f(x)={\widehat E}_x[f({\widehat \xi}^+(t\wedge{\widehat \sigma}_0))],
\end{equation}
where $\sigma_0$ and ${\widehat \sigma}_0$ are the first hitting time to $0$ of
 $\xi^+(t)$ and ${\widehat \xi}^+(t)$, respectively.
Introduce, further,  the semigroups for {\em processes with  extinction on hitting  $0$} by
\begin{equation}
T^0_t f(x)=E_x\left[f(\xi^+(t))\cdot{\bf 1}_{[t<\sigma_0]}\right] \ \mbox{and} \ \  {\widehat T}^0_t f(x)={\widehat E}_x\left[f({\widehat \xi}^+(t))\cdot{\bf 1}_{[t<{\widehat \sigma}_0]}\right].
\end{equation}
$T^-_t$ and ${\widehat T}^-_t$ are Markovian semigroups and $T^0_t$ and ${\widehat T}^0_t$ are sub-Markovian semigroups.  Note also that $T^+_t$ , ${\widehat T}^+_t$, $T^0_t$ and ${\widehat T}^0_t$ have the  strong Feller property but $T^-_t$ and ${\widehat T}^-_t$ have the Feller property only.  It holds that
\begin{equation}
T^0_t f=T^-_t({\bf 1}_{(0,\infty)}\cdot f) \quad \mbox{and}\quad {\widehat T}^0_t f={\widehat T}^-_t({\bf 1}_{(0,\infty)}\cdot f).
\end{equation}

We have the following duality relations which form another key lemma in the proof of Th.1.2:
\begin{lm} \ \ For $x, y\in [0,\infty)$ and $t>0$,
\begin{equation}
T^+_t{\bf 1}_{[0,y]}(x)={\widehat T}^0_t{\bf 1}_{[x,\infty)}(y) \quad\mbox{and}\quad T^-_t{\bf 1}_{[0,y]}(x)={\widehat T}^+_t{\bf 1}_{[x,\infty)}(y).
\end{equation}
More generally, for $x, y\in [0,\infty)$ and $0\leq t_0<t_1<\ldots<t_{2n-1}<t_{2n}<t_{2n+1}$,
\begin{eqnarray}
&  &T^+_{t_1-t_0}T^-_{t_2-t_1}T^+_{t_3-t_2}\cdots T^+_{t_{2n-1}-t_{2n-2}}T^-_{t_{2n}-t_{2n-1}}{\bf 1}_{[0,y]}(x)\\
&=&{\widehat T}^+_{t_{2n}-t_{2n-1}}{\widehat T}^0_{t_{2n-1}-t_{2n-2}}\cdots {\widehat T}^0_{t_3-t_2}{\widehat T}^+_{t_2-t_1}{\widehat T}^0_{t_1-t_0}{\bf 1}_{[x,\infty)}(y), \nonumber
\end{eqnarray}
and
\begin{eqnarray}
&  &T^+_{t_1-t_0}T^-_{t_2-t_1}T^+_{t_3-t_2}\cdots T^+_{t_{2n-1}-t_{2n-2}}T^-_{t_{2n}-t_{2n-1}}T^+_{t_{2n+1}-t_{2n}}{\bf 1}_{[0,y]}(x)\\
&=&{\widehat T}^0_{t_{2n+1}-t_{2n}}{\widehat T}^+_{t_{2n}-t_{2n-1}}{\widehat T}^0_{t_{2n-1}-t_{2n-2}}\cdots {\widehat T}^0_{t_3-t_2}{\widehat T}^+_{t_2-t_1}{\widehat T}^0_{t_1-t_0}{\bf 1}_{[x,\infty)}(y).\nonumber
\end{eqnarray}
\end{lm}

Admitting this lemma for a moment, we now proceed to prove Th. 1.2.\medskip

\noindent {\em Proof of Th. 1.2.}  \ \ Let $F=[t_0,t_1]\cup [t_2,t_3]\ldots \cup[t_{2n-2},t_{2n-1}]$ be an elementary set in $[0,1]$ and $({\bf X}, {\bf X}')$ be a $(1^-,F)$-coupling of the Harris flow.  Set $\xi(t)=X_{0,t}(0)-X'_{0,t}(0)$.  Then $|\xi(t)|$ is a time-inhomogeneous diffusion process which behaves as a reflecting $L$-diffusion when $t\in F$ and as an absorbing $L$-diffusion (i.e., $L$-diffusion with $0$ as a trap) when $t\in [0,1]\setminus F$.  It is known that $P(S_X\ni t)=0$ for every $t\in [0,1]$ (cf. [T 2]). Then (2.8), combined with this remark, yields that
$$P(|S_X\cap F|=\infty)=P(S^{acc}_X\cap F\ne \emptyset)=\frac12 E[|\xi(1)|^2]. $$By applying the It\^o formula for $\xi(t)$ on each interval $[t_k,t_{k+1}]$, we have 
$$\frac12 E[|\xi(1)|^2]=\int^1_0 E[(1-b)(\xi(t))]dt=1-\int^1_0 E[b(\xi(t))]dt,$$and hence,
\begin{equation}
P(S^{acc}_X\cap F= \emptyset)=\int^1_0 E[b(\xi(t))]dt.
\end{equation}
On the other hand, 
$$E[b(\xi(t))]=\left\{\begin{array}{ll}T^+_{t_1-t_0}T^-_{t_2-t_1}\cdots T^-_{t_{2k}-t_{2k-1}}T^+_{t-t_{2k}}b(0), & \ \mbox{if} \ \ t_{2k}\leq t<t_{2k+1}\\
T^+_{t_1-t_0}T^-_{t_2-t_1}\cdots T^+_{t_{2k-1}-t_{2k-2}}T^-_{t-t_{2k-1}}b(0), &  \ \mbox{if} \ \ t_{2k-1}\leq t<t_{2k}\end{array}\right..$$
Noting  $b(x)=\int_{[0,\infty)}{\bf 1}_{[0,y]}(x)\mu(dy)$, we have by Lemma 3.1 the following:
\begin{eqnarray*}& &E[b(\xi(t))]\\
&=&\left\{\begin{array}{ll}\int_0^\infty \mu(dy)({\widehat T}^0_{t-t_{2k}}{\widehat T}^+_{t_{2k}-t_{2k-1}}\cdots {\widehat T}^+_{t_2-t_1}{\widehat T}^0_{t_1-t_0}{\bf 1}_{[0,\infty)})(y), &  \mbox{if} \ \ t_{2k}\leq t<t_{2k+1}\\
\int_0^\infty \mu(dy)({\widehat T}^+_{t-t_{2k-1}}{\widehat T}^0_{t_{2k-1}-t_{2k-2}}\cdots {\widehat T}^+_{t_2-t_1}{\widehat T}^0_{t_1-t_0}{\bf 1}_{[0,\infty)})(y), &   \mbox{if} \ \ t_{2k-1}\leq t<t_{2k}\end{array}\right.
\end{eqnarray*}
If the random set ${\widetilde S}$ is defined by (1.7), it is not difficult to deduce, from the last expression of $E[b(\xi(t))]$, that $\int^1_0 E[b(\xi(t))]dt$ coincides with $P({\widetilde S}\cap F=\emptyset)$.  Then  $P({\widetilde S}\cap F=\emptyset)=P(S^{acc}_X\cap F= \emptyset)$ by (3.8). Since this holds for every elementary set $F$, we can conclude that $S^{acc}_X\stackrel{d}{=}{\widetilde S}$.
\hfill \qed \bigskip

\noindent{\em Proof of Lemma 3.1.}  \ \ First, we prove (3.5). For  $\lambda>0$, let $U^+_\lambda$ and $\hat{U}_\lambda^0$ be the resolvent operators associated with the semigroups $T^+_t$ and  $\hat{T}_t^0$ respectively. Let $f$ be continuous and compactly supported in $(0,\infty)$. Then $u=U^+_\lambda f$ solves Poisson's equation  
$$  Lu-\lambda u=-f,$$
with the boundary conditions $u^\prime(0+)=u(\infty)=0$. Define functions $g$ and $v$ via   
$$g(y)= \int_0^y \frac{f(x)}{a(x)}dx \quad \mbox{ and } \quad v(y)= \int_0^y \frac{u(x)}{a(x)}dx, $$
where $a(x)=(1-b(x))$.  Dividing Poisson's equation through by $a(x)$ and integrating, we obtain 
$$\hat{L}v-\lambda v=-g.$$
Moreover $v$ and $g$ are bounded and $v(0)=0$. Thus we must have $v=\hat{U}_\lambda^0 g.$  Letting $f$ approach a delta function we may write the relationship between $u$ and $v$ as:
$$ \frac{1}{a(z)}\hat{U}_\lambda^0{\bf 1}_{[z,\infty)}(y)
=\int_0^y \frac{u^+_\lambda(x,z)}{a(x)}dx, $$
where $u^+_\lambda$ is the continuous version of the  resolvent density corresponding to $U^+_\lambda$. Recalling the symmetry relation,
$$\frac{1}{a(x)}u^+_\lambda(x,z)a(z)=u^+_\lambda(z,x),$$
we obtain
$$ \hat{U}_\lambda^0{\bf 1}_{[z,\infty)}(y)=\hat{U}_\lambda^+{\bf 1}_{[0,y]}(z),$$
from which the first equality of (3.5) follows by uniqueness of Laplace transforms. The second equality may be proved by a similar method.

(3.6) and (3.7) can be proved by applying (3.5) successively: For example,
\begin{eqnarray*}
&  & T^+_{t_1-t_0}\cdot T^-_{t_2-t_1}{\bf 1}_{[0,y]}(x)=\int_{[0,\infty)}T^+_{t_1-t_0}(x,du) T^-_{t_2-t_1}{\bf 1}_{[0,y]}(u)\\
&=&\int_{[0,\infty)}T^+_{t_1-t_0}(x,du){\widehat T}^+_{t_2-t_1}{\bf 1}_{[u,\infty)}(y)=\int\!\!\!\int_{0\leq u\leq v<\infty}T^+_{t_1-t_0}(x,du){\widehat T}^+_{t_2-t_1}(y,dv)\\
&=&\int_{[0,\infty)}{\widehat T}^+_{t_2-t_1}(y,dv)T^+_{t_1-t_0}{\bf 1}_{[0,v]}(x)=\int_{[0,\infty)}{\widehat T}^+_{t_2-t_1}(y,dv){\widehat T}^0_{t_1-t_0}{\bf 1}_{[x,\infty)}(v)\\
&=&{\widehat T}^+_{t_2-t_1}\cdot {\widehat T}^0_{t_1-t_0}{\bf 1}_{[x,\infty)}(y).
\end{eqnarray*}
This proves a particular case of (3.6). In the same way, the general case can be proved easily by induction.  \hfill \qed
\begin{rem}
We remark that an alternative proof of (3.5) is possible by means of the time reversal of stochastic flows on the half line. A stochastic flow on the half line $[0,\infty)$ is defined similarly by replacing the whole line ${\bf R}$ by $[0,\infty)$ in Def.1.2.  A key idea in the proof is to construct a stochastic flow ${\bf X}=(X_{s,t})$ on $[0,\infty)$ whose one-point motion $t\mapsto X_{0,t}(x), x\in{\bf R}$, is given by the absorbing $L$-diffusion $\xi^-(t)$, i.e., the diffusion with the semigroup $T_t^-$, and then show that its {\em time reversed flow} ${\widehat {\bf X}}=({\widehat X}_{s,t})$, defined by ${\widehat X}_{s,t}=\left(X_{-t,-s}\right)^{-1}$, has the one-point motion given by the reflecting ${\widehat L}$-diffusion ${\widehat \xi}^+(t)$, i.e., the diffusion with the semigroup ${\widehat T}_t^+$. Here, for a right-continuous and non-decreasing $\varphi:[0,\infty)\to[0,\infty)$ such that $\lim_{x\nearrow\infty}\varphi(x)=\infty$, $\varphi^{-1}$ is the right-continuous inverse of $\varphi$: $\varphi^{-1}(x)=\inf\{y|\varphi(y)>x\}$.  This is connected to the fact that $L$ and ${\widehat L}$, when written in H\"ormander form, differ only in the sign of the drift term.  The corresponding fact in the case of stochastic flows of homeomorphisms is well-known (cf. [K] p.131, [IW] p.265).
\end{rem}
\setcounter{section}{3}
\setcounter{equation}{0}
\section{Proof of Th. 1.3.}  \ \ Consider a Harris flow ${\bf X}$ satisfying (1.4), (1.6) and (1.11). \medskip

\noindent{\em Proof of Cor. 1.1.} \ \ It is sufficient to show that the set of zeros of ${\widehat L}$-diffusion ${\widehat \xi}(t)$ has the Hausdorff dimension $(1-\a)/(2-\a)$, ${\widehat P}_0$-almosy surely.  The  set of zeros of ${\widehat \xi}(t)$ is the range of the inverse local time $l^{-1}(t)$ at $0$ of ${\widehat \xi}(t)$, which is a subordinator with exponent $\Psi(\lambda)=g_\lambda(0,0)^{-1}$:
$$  E\left(e^{-\lambda l^{-1}(t)}\right)=e^{-t\Psi(\lambda)}=e^{-t/g_\lambda(0,0)}.  $$
Here, $g_\lambda(x,y)$  is the Green function (resolvent density) with respect to the speed measure $dx$ of reflecting ${\widehat L}$-diffusion where ${\widehat L}=\frac{d}{dx}(1-b(x))\frac{d}{dx}$.  If we introduce the scale $\xi=\int_0^x(1-b(y))^{-1}dy$  as the coordinate of $[0,\infty)$, then ${\widehat L}=(1-{\tilde b}(\xi))^{-1}\frac{d^2}{d^2\xi}$  where ${\tilde b}(\xi)=b(x(\xi))$, so that the speed measure in the new coordinate is given by $d{\tilde m}(\xi)=a(\xi)d\xi$  with  $a(\xi)=1-{\tilde b}(\xi)$.   It is easy to deduce from (1.11) that $a(\xi)\asymp \xi^{\a/(1-\a)}$  as $\xi\to 0$.  Let ${\tilde g}_\lambda(\xi,\eta)$ be the Green function for  ${\widehat L}$-diffusion with respect to the speed measure so that ${\tilde g}_\lambda(0,0)=g_\lambda(0,0)$.   By Th.2.3 in p.243 of [KW], we have
$$\Psi(\lambda)={\tilde g}_\lambda(0,0)^{-1}\asymp\lambda^{1/(2+\frac{\a}{1-\a})}=\lambda^{\frac{1-\a}{2-\a}}\quad\mbox{as} \ \ \lambda\to \infty.$$
Then we can conclude that the range of the subordinator  $l^{-1}(t)$ has the Hausdorff dimension $\frac{1-\a}{2-\a}$ almost surely, by a result of Blumenthal and Getoor (cf. [B], p. 94, Th. 16).  \hfill \qed \medskip

Now we proceed to prove Th. 1.3. We need several lemmas.
\begin{lm} (i) \   Let $\Phi_1, \Phi_2 \in L_2^{us}({\cal F}^X_{-\infty,\infty})$ and consider their linear combination $\Phi=\a\Phi_1+\b\Phi_2\in L_2^{us}({\cal F}^X_{-\infty,\infty})$.  If $A\in {\cal B}({\cal C})$ satisfies  
$P(S_{\Phi_1}\in A)=P(S_{\Phi_2}\in A)=1$,
then it holds that
$P(S_{\Phi}\in A)=1.$

\noindent (ii) \   Let $\Phi_n \in L_2^{us}({\cal F}^X_{-\infty,\infty})$, $n=1,2,\ldots,$ constitute a dense family in $L_2^{us}({\cal F}^X_{-\infty,\infty})$.  If $A\in {\cal B}({\cal C})$ satisfies  
$P(S_{\Phi_n}\in A)=1 \quad \mbox{for all} \ \ n$,
then it holds that
$P(S_{\Phi}\in A)=1 \quad \mbox{for all} \ \ \Phi \in L_2^{us}({\cal F}^X_{-\infty,\infty}).$
\end{lm}
{\em Proof.} \ \ 
According to Theorem 3d12 of [T 5], every $A\in {\cal B}({\cal C})$ is associated with a closed subspace ${\cal H}_A$ 
of $L_2({\cal F}^X_{-\infty,\infty})$  such that the spectral measure $\mu_\Phi$ of any $\Phi$ satisfies
\[
|| P_A \Phi ||^2= \mu_\Phi(A),\]
where $P_A$ denotes the orthogonal projection onto ${\cal H}_A$.  Both parts of this lemma are immediate consequences.\hfill \qed
\begin{lm} \ \  Let $t_1<t_2<t_3$ and $\Phi=\Phi_1\Phi_2\in L_2^{us}({\cal F}^X_{t_1,t_3})$ such that $\Phi_1\in L_2^{us}({\cal F}^X_{t_1,t_2})$ and $\Phi_2\in L_2^{us}({\cal F}^X_{t_2,t_3})$. Then,
$$ S_\Phi\cap [t_1,t_2]\stackrel{d}{=}S_{\Phi_1}, \quad S_\Phi\cap [t_2,t_3]\stackrel{d}{=}S_{\Phi_2}.$$
Furthermore, $S_\Phi\cap [t_1,t_2]$ and $S_\Phi\cap [t_2,t_3]$ are mutually independent.
\end{lm}
The proof is easy and omitted.
\begin{lm}
Let $S$ be a ${\cal C}_{[0,1]}$-valued random variable and assume, for $0<\b<1$ and $K>0$, that
$$ P(S\cap [t, t+\epsilon]\ne \emptyset)\leq K\epsilon^\b\quad \mbox{for all} \ \ 0<\epsilon<1\quad \mbox{and} \ \ t\in[0,1].$$
Then, $P(\dim S\leq 1-\b)=1$.
\end{lm}
{\em Proof}. \ \ For every $a>1-\b$, we have
\begin{eqnarray*}
& &E\left(\sum_{k=1}^n {\bf 1}_{\{S\cap [\frac{k-1}{n}, \frac{k}{n}]\ne \emptyset\}}\cdot \left(\frac{1}{n}\right)^a\right)=\sum_{k=1}^n  P\left(S\cap \left[\frac{k-1}{n}, \frac{k}{n}\right]\ne \emptyset\right)\cdot\left(\frac{1}{n}\right)^a\\
&\leq& nK\left(\frac{1}{n}\right)^\b\cdot \left(\frac{1}{n}\right)^a=K\cdot n^{1-(\b+a)} \to 0 \quad\mbox{as} \quad n\to \infty.
\end{eqnarray*}
Hence, there exists a subsequence $n_\nu\to \infty$ such that, almost surely, 
$$ \sum_{k=1}^{n_\nu} {\bf 1}_{\{S\cap [\frac{k-1}{n_\nu}, \frac{k}{n_\nu}]\ne \emptyset\}}\cdot \left(\frac{1}{n_\nu}\right)^a\to 0 \quad\mbox{as}\quad \nu\to \infty.$$
Let ${\cal C}_\nu$ be the collection of  intervals $E_k=[\frac{k-1}{n_\nu}, \frac{k}{n_\nu}]$, \ $k=1,\ldots,n_\nu,$  \ which have nonempty intersections with  the set $S$. Then ${\cal C}_\nu$ is a covering of $S$ and
$$ \sum_{E_k\in {\cal C}_\nu}({\rm diam} \ E_k)^a \to 0 \quad\mbox{a.s., \ as}\quad \nu\to \infty.$$   Hence, $\dim S\leq a$, a.s., implying that $\dim S\leq 1-\b$, a.s.
\hfill \qed \medskip

\noindent{\em Proof of Th. 1.3}.  \ \ It is sufficient to show that
\begin{equation}
\dim S_{\Phi}\leq \frac{1-\a}{2-\a}\quad\mbox{a.s.}
\end{equation}
for $\Phi\in L_2^{us}({\cal F}_{0,1}^X)$. Indeed, if (4.1) is true for $\Phi\in L_2^{us}({\cal F}_{0,1}^X)$, then by the stationarity of the flow, it is also true for $\Phi\in L_2^{us}({\cal F}_{n,n+1}^X)$. By Lemma 4.2, (4.1) is true for a finite product of such $\Phi$'s.  Since linear combinations of such products are dense in $L_2({\cal F}_{-\infty,\infty}^X)$, we can conclude by Lemma 4.1 that (4.1) is true for any $\Phi\in L_2^{us}({\cal F}_{-\infty,\infty}^X)$.

First, we consider the case when $\Phi\in L_2^{us}({\cal F}_{0,1}^X)$ is given by
$$\Phi=f(X_{0,1}(x_1),\ldots,X_{0,1}(x_n)), \quad x_1,\ldots, x_n\in {\bf R},$$ and a function $f$ is uniformly Lipschitz-continuous on ${\bf R}^n$.  

Let $F=[t,t+\epsilon]$, $0\leq t<t+\epsilon\leq 1$, and let $({\bf X}, {\bf X}')$ be a $(1^-, F)$-joining.  Then we know by Lemma 2.2 that  $2P(S^{acc}_X\cap F\ne\emptyset)=E(|X_{0,1}(0)-X'_{0,1}(0)|^2)$ 
and similarly, we have 
$2P(S^{acc}_\Phi\cap F\ne\emptyset)=E(|\Phi-\Phi'|^2)$ 
where $\Phi'=f(X'_{0,1}(x_1),\ldots,X'_{0,1}(x_n))$. Therefore, noting that $E(|X_{0,1}(x)-X'_{0,1}(x)|^2)$ is independent of $x$, we have
\begin{eqnarray}
& &P(S^{acc}_\Phi\cap F\ne\emptyset)=\frac12E(|\Phi-\Phi'|^2)\nonumber\\
&\leq& KE(|X_{0,1}(0)-X'_{0,1}(0)|^2)=2KP(S^{acc}_X\cap F\ne\emptyset)
\end{eqnarray}
where a constant $K$ depends on $n$ and the Lipschitz constant of $f$.

Let $\{{\widehat \xi}^+(t), {\widehat P}_\xi\}$  be the reflecting ${\widehat L}$-diffusion on $[0,\infty)$. As in the proof of Cor.1.1, {\em take a canonical scale $\xi$ as the coordinate} so that  ${\widehat L}=\frac{d^2}{a(\xi)d\xi^2}$ and we have
$ a(\xi)\asymp \xi^{\a/(1-\a)}$ as $\xi \to 0$ and $a(\xi)\to 1$ as $\xi \to \infty$.  Let $\mu(d\xi)=da(\xi)$.  By what we have shown above,
\begin{eqnarray*}& &P(S^{acc}_X\cap [t,t+\epsilon]\ne\emptyset)=P({\widetilde S}\cap [t,t+\epsilon]\ne \emptyset)\\
&=&\int_0^1 {\widehat P}_\mu\left({\widehat \xi}^+(u-s)=0 \ \mbox{for some} \ s\in[0,u]\cap [t,t+\epsilon]\right)du\\
&=&\int_t^1 {\widehat P}_\mu\left({\widehat \xi}^+(\theta)=0 \ \mbox{for some} \ \theta\in[(u-t-\epsilon)_+,u-t]\right)du\\
&=&O(\epsilon)+\int_t^1 {\widehat P}_\mu\left({\widehat \xi}^+(\theta)=0 \ \mbox{for some} \ \theta\in [u-t,u-t+\epsilon]\right)du.
\end{eqnarray*}
We would show  
\begin{equation}I(t):=\int_t^1 {\widehat P}_\mu\left({\widehat \xi}^+(\theta)=0 \ \mbox{for some} \ \theta\in [u-t,u-t+\epsilon]\right)du=O\left(\epsilon^{\frac{1}{2-\a}}\right)\end{equation}
 as $\epsilon\to 0$ uniformly in $t\in[0,1].$  If we can show this, then
$$P(S^{acc}_X\cap [t,t+\epsilon]\ne\emptyset)=O(\epsilon^{1/(2-\a)})$$  as $\epsilon\to 0$ uniformly in $t\in[0,1]$ and, combining this with (4.2), we see that $P(S^{acc}_\Phi\cap [t,t+\epsilon]\ne\emptyset)=O(\epsilon^{1/(2-\a)})$, so that, by Lemma 4.3, we can conclude that the estimate (4.1) holds for $\Phi$  because $1-1/(2-\a)=(1-\a)/(2-\a)$.

To obtain (4.3), we estimate 
\begin{eqnarray*}
 & & I(t)\leq \int_0^1 {\widehat P}_\mu\left({\widehat \xi}^+(\theta)=0 \ \mbox{for some} \ \theta\in [u,u+\epsilon]\right)du\\
 &=&\int_0^1{\widehat E}_\mu\left({\widehat P}_{{\widehat \xi}^+(u)}[{\widehat \sigma}_0\leq \epsilon]\right)du\leq e\int_0^1e^{-u}{\widehat E}_\mu\left({\widehat P}_{{\widehat \xi}^+(u)}[{\widehat \sigma}_0\leq \epsilon]\right)du\\
&\leq& e\int_0^\infty e^{-u}{\widehat E}_\mu\left({\widehat P}_{{\widehat \xi}^+(u)}[{\widehat \sigma}_0\leq \epsilon]\right)du= e\int_{[0,\infty)}\mu(d\xi)\int_{[0,\infty)}{\tilde g}_1(\xi,\eta){\widehat P}_\eta[{\widehat \sigma}_0\leq \epsilon]a(\eta)d\eta,
\end{eqnarray*}
where ${\widehat \sigma}_0$ is the first hitting time of ${\widehat \xi}^+(t)$ to $0$. Since the resolvent density ${\tilde g}_1(\xi,\eta)$ is bounded, we have, for some $C>0$,
$$ I(t)\leq C\int_{[0,\infty)}{\widehat P}_\eta[{\widehat \sigma}_0\leq \epsilon]a(\eta)d\eta.$$
The process  ${\widehat \xi}^+(t)$ under ${\widehat P}_\eta$,  $\eta>0$, and in the coordinate $\xi$, is obtained from a one-dimensional Brownian motion $B(t)$ with $B(0)=0$ by
$${\widehat \xi}^+(t)=|\eta+B(A^{-1}(t))|\quad \mbox{where} \ \ A(t)=\int_0^t a(|\eta+B(s)|)ds.$$
Hence, 
$${\widehat P}_\eta({\widehat \sigma}_0\leq \epsilon)=P\left(\int_0^{\sigma_0} a(|\eta+B(s)|)ds\leq \epsilon\right) \quad \mbox{where} \ \ \sigma_0=\min\{s|\eta+B(s)=0\}.$$
and, noting $a(\xi)\geq K^{-1}\cdot\xi^{\a/(1-\a)}\wedge 1$ for some $K>0$,
$${\widehat P}_\eta({\widehat \sigma}_0\leq \epsilon)\leq P\left(\int_0^{\sigma_0} \left(|\eta+B(s)|^{\a/(1-\a)}\wedge 1\right) ds\leq K\epsilon\right).$$
The scaling property of $B(t)$ combined with an easy inequality $(\epsilon a)\wedge1\geq \epsilon(a\wedge 1)$ for $a>0$ and $1\geq \epsilon >0$ yields that the RHS is dominated by $\phi(\epsilon^{-(1-\a)/(2-\a)}\eta)$, where
$$\phi(\eta)= P\left(\int_0^{\sigma_0} \left(|\eta+B(s)|^{\a/(1-\a)}\wedge 1\right) ds\leq K\right).$$
Then, 
\begin{eqnarray*}& &I(t)\leq C\int_{[0,\infty)}\phi(\epsilon^{-(1-\a)/(2-\a)}\eta)a(\eta)d\eta\\
&\leq& K'\int_{[0,\infty)}\phi(\epsilon^{-(1-\a)/(2-\a)}\eta)\eta^{\a/(1-\a)}d\eta= K'\epsilon^{1/(2-\a)}\int_{[0,\infty)}\phi(\eta)\eta^{\a/(1-\a)}d\eta
\end{eqnarray*}
and we have otained (4.3).\medskip

In the same way, we have the estimate (4.1) for $\Phi=f(X_{s,t}(x_1),\ldots,X_{s,t}(x_n))$, $x_1,\ldots, x_n\in {\bf R}$, $s<t$, where $f$ is uniformly Lipschitz continuous on ${\bf R}^n$. Then, by Lemma 4.2, we have the estimate (4.1) for $\Phi=\Phi_1 \Phi_2\cdots\Phi_m \in L_2^{us}({\cal F}_{0,1}^X)$ if $t_0=0<t_1<t_2<\cdots<t_m=1$, and $\Phi_k\in {\rm ub}[L^2({\cal F}_{t_{k-1},t_k}^X)]$, $k=1,2,\ldots,m$, is given in the form \[
\Phi_k=f_k\left(X_{t_{k-1},t_k}(x_1^{(k)}),\ldots,X_{t_{k-1},t_k}(x_{n_k}^{(k)})\right),
\]
 $x_1^{(k)},\ldots, x_{n_k}^{(k)}\in {\bf R}$, where $f_k$ is uniformly Lipschitz continuous on ${\bf R}^{n_k}$. By Lemma 4.1 (i), the estimate (4.1) still holds for a finite linear combination of such functionals and this class of functionals is dense in $L_2^{us}({\cal F}_{0,1}^X)$.\hfill \qed

\begin{flushleft}
Jon Warren\\
Department of Statistics, University of Warwick\\
Coventry, CV4 7AL, United Kingdom
\end{flushleft}

\begin{flushleft}
Shinzo Watanabe\\
527-10, Chaya-cho, Higashiyama-ku\\
Kyoto, 605-0931, Japan
\end{flushleft}

\end{document}